\documentclass[12pt]{article}
\title{\sc more indecomposable polyhedra}
\usepackage{ graphicx}
\usepackage{color}
\usepackage{hyperref}
\usepackage[adobe-utopia]{mathdesign}
\usepackage[T1]{fontenc}
\let\svthefootnote\thefootnote
\newcommand\blankfootnote[1]{ \let\thefootnote\relax\footnotetext{#1}
 \let\thefootnote\svthefootnote
}
\date{}
\author{\sc Krzysztof Przes\l awski \& David Yost}
\begin{document}
\maketitle
\newtheorem{pr}{\sc Proposition}
\newtheorem{thm}[pr]{\sc Theorem}
\newtheorem{lem}[pr]{\sc Lemma}
\newtheorem{cor}[pr]{\sc Corollary}
\def\pf{\nt{\bf Proof.\ }}
\def\al{\alpha} \def\be{\beta}
\def\de{\delta}
\def\la{\lambda}
\def\centreline{\centerline}
\def\h{{1\over2}}
\def\co{{\rm co}}

\def\qed{\hfill $\diamondsuit$}
\def\er{\mathbb{R}}
\def\erde{\er^d}
\def\rthree{\er^3}
\def\proof{\noindent{\bf Proof.\/}\ }
\def\pf{\proof}

\centreline{\it Dedicated to the memory of Carlos Ben{\'\i}tez}

\bigskip

\bigskip

\begin{abstract}
 We apply combinatorial methods to a geometric problem: the classification of polytopes, in terms of Minkowski decomposability. Various properties of skeletons of polytopes are exhibited, each sufficient to guarantee indecomposability of a significant class of polytopes. We illustrate further the power of these techniques, compared with the traditional method of examining triangular faces,  with several applications. In any dimension $d\ne2$, we show that of all the polytopes with $d^2+\h d$ or fewer edges, only one is decomposable. In 3 dimensions, we complete the classification, in terms of  decomposability, of the 260 combinatorial types of polyhedra with 15 or fewer edges.
\end{abstract}
\bigskip
\blankfootnote{Mathematics subject classification: 52B10, 52B11, 52B05; Keywords: polytope, decomposable}

What happens if you have two line segments in the plane, oriented in different directions, and you calculate all the sums of all pairs of elements, one from each segment? Of course you end up with a rectangle, or at least a parallelogram. Do the same again with a triangle and a line segment in three dimensions: this time, you get a  prism. Thus the prism and the parallelogram are {\it decomposable}; they can be expressed as the (Minkowski) sum of two dis-similar convex bodies. (Recall that that two polytopes
are similar if one can be obtained from the other by a dilation and a translation.) On the other hand, any triangle, tetrahedron or octahedron is {\it indecomposable}.

We refer to \cite{Gr} for a
general  introduction to the theory of polytopes, as well as for specific results.
Determining the decomposability of a polytope can be reduced to a computational problem in linear algebra \cite{M,S}. That is,
given the co-ordinates of its vertices, all we have to do is calculate the rank of a rather large matrix.  However that is not the approach to be taken here.

The  edges and vertices  of any polytope obviously constitute a graph, sometimes known as its skeleton. In the case of a polyhedron, this will be isomorphic to a planar graph.
All the geometric conclusions of this paper will be established by considering the properties of this graph. Section 1 develops a number of sufficient conditions for indecomposability (or decomposability).
Our results have wider applicability than earlier results in this area, which generally relied on the existence large families of triangular faces. Section 2 applies them to complete the classification of 3-dimensional polyhedra with up to 15 edges. Section 3 applies them to completely classify, as indecomposable or decomposable, all $d$-dimensional polytopes with up to $d^2+\h d$ edges. We also show that there is no $d$-dimensional polytope at all with  $2d$ vertices and $d^2+1$ edges, for $d\ne3$.

\section{Geometric graphs and indecomposability}

We will not give a thorough history of this topic, but it is important to recall
some preliminary information.

We depend heavily on the concept of a {\it geometric graph}, which was pioneered by
Kallay \cite{K}. He defined a geometric graph as any graph $G$ whose vertex set $V$
is a subset of a finite-dimensional real vector space $X$, and whose edge set $E$ is a subset of  the line segments
joining members of $V$. (Of course, $X$ will be isomorphic to $\erde$ for some $d$, but we prefer this basis-free formulation.) It is largely a formality whether we consider an edge to be an unordered pair or a line segment.
It is significant that such a graph  need not be the edge graph of any polytope. He then extended the notion of decomposability to such
graphs in the following manner.

For convenience, let us say that a function $f\colon V\to X$ is a {\it decomposing\/  function} for the graph $(V,E)$
if it has the property that $f(v)-f(w)$ is a scalar multiple of $v-w$ for each edge $[v,w]\in E$. (This is slightly different from Kallay's {\it local similarity}; he insisted on strictly positive multiples.)
A geometric graph $G=(V,E)$ is then called {\it decomposable} if there is a decomposing function which is neither constant, nor the restriction of a homothety on {$X$}. If the only non-constant decomposing functions are homotheties then $G$ is called {\it indecomposable}.

Significantly, Kallay showed \cite[Theorem 1]{K} that a polytope is indecomposable if and only if
its edge graph is indecomposable in this sense. Exploiting an idea of McMullen \cite{McM} and Kallay \cite[Theorem 1b]{K}, we showed in \cite[Theorem 8]{PY} that it is even sufficient just to have an indecomposable subgraph which contains at least one vertex from every facet (maximal face).
A strategy for proving indecomposability of a polytope is thus to prove that certain simple geometric graphs are indecomposable, and by building up to show that the entire skeleton of our polytope is indecomposable. (It also would be interesting to formulate somehow a notion of {\it primitivity} for such graphs.)

Building on the concept introduced in \cite[p. 139]{Yo}, let us say that a geometric graph $G=(V,E)$ is a simple extension of a  geometric graph  $G_0=(V_0,E_0)$ if $G$ has one more vertex and two more edges than $G_0$. More precisely, we mean that there is a unique $v\in V\setminus V_0$, and distinct vertices $u$ and $w$ in $V_0$, such  that $E=E_0\cup \{[u,v]\}\cup\{[v,w]\}$. Observe that the existence of these two edges means that the value of any decomposing function at $v$ is determined by its values at $u$ and $w$. No assumption is made about whether $[u,w]$ is an edge of either graph. Our first result is a special case of the next one, but it is so useful and so easy to prove that it is worth stating separately.

\begin{pr}\label{extend}
Suppose that $G_0,G_1,\ldots, G_n$ are  geometric graphs, that $G_{i+1}$ is a simple extension of $G_{i}$ for each $i$,  and that $G_0$ is indecomposable.
 Then $G_n$ is also indecomposable.
\end{pr}

\proof  It is clearly sufficient to prove this when $n=1$, and this follows from the observation in the preceding paragraph. \qed

\medskip
Let us illustrate how this  can be applied in the simplest cases, polyhedra for which ``sufficiently many"  faces are triangles \cite[\S3]{Sh}.  Any edge is obviously indecomposable, and then Proposition \ref{extend} easily implies that any triangle is indecomposable. Furthermore if an indecomposable geometric graph shares an edge with a triangle, then their union is easily proved to be indecomposable. It follows that the union of a
  {\it  chain} of triangles, as defined in \cite[p. 92]{Sh}, is an indecomposable graph.
This makes it clear that a polyhedron must be indecomposable if every face is a triangle. If every face but one is a triangle, it remains true that the triangular faces can be ordered into a  chain, whose union is the entire skeleton of the polyhedron;  again indecomposability is assured. The same holds if all faces but two are  triangular, and the non-triangular faces do not share an edge. If all faces but two are  triangular, but the non-triangular faces do share an edge, then the triangular faces can still be ordered into a  chain, whose union will contain every vertex of the polyhedron and every edge but one. So, indecomposability is assured, whenever there are two or fewer non-triangular faces.

This conclusion no longer holds if we have three non-triangular faces, as the triangular prism is  decomposable. On the other hand, there are also many indecomposable polyhedra with precisely three non-triangular faces. For comparison,
let us mention that that a polyhedron with only three, or fewer, triangular faces is automatically decomposable  \cite[\S 6]{S}.

\begin{figure}
\hspace{0.07\textwidth}
\includegraphics[width=0.4\textwidth]{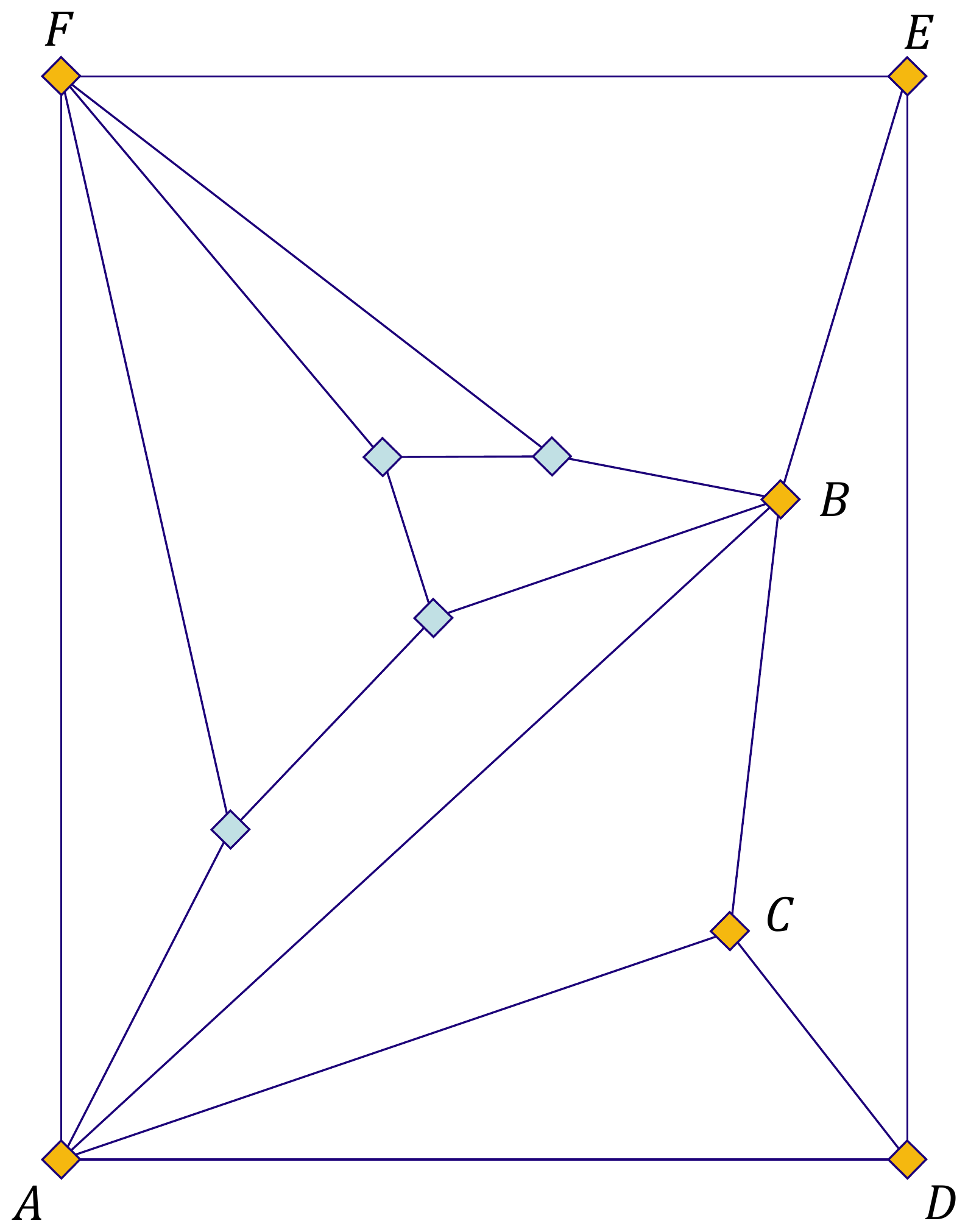}
\hspace{0.06\textwidth}
\includegraphics[width=0.4\textwidth]{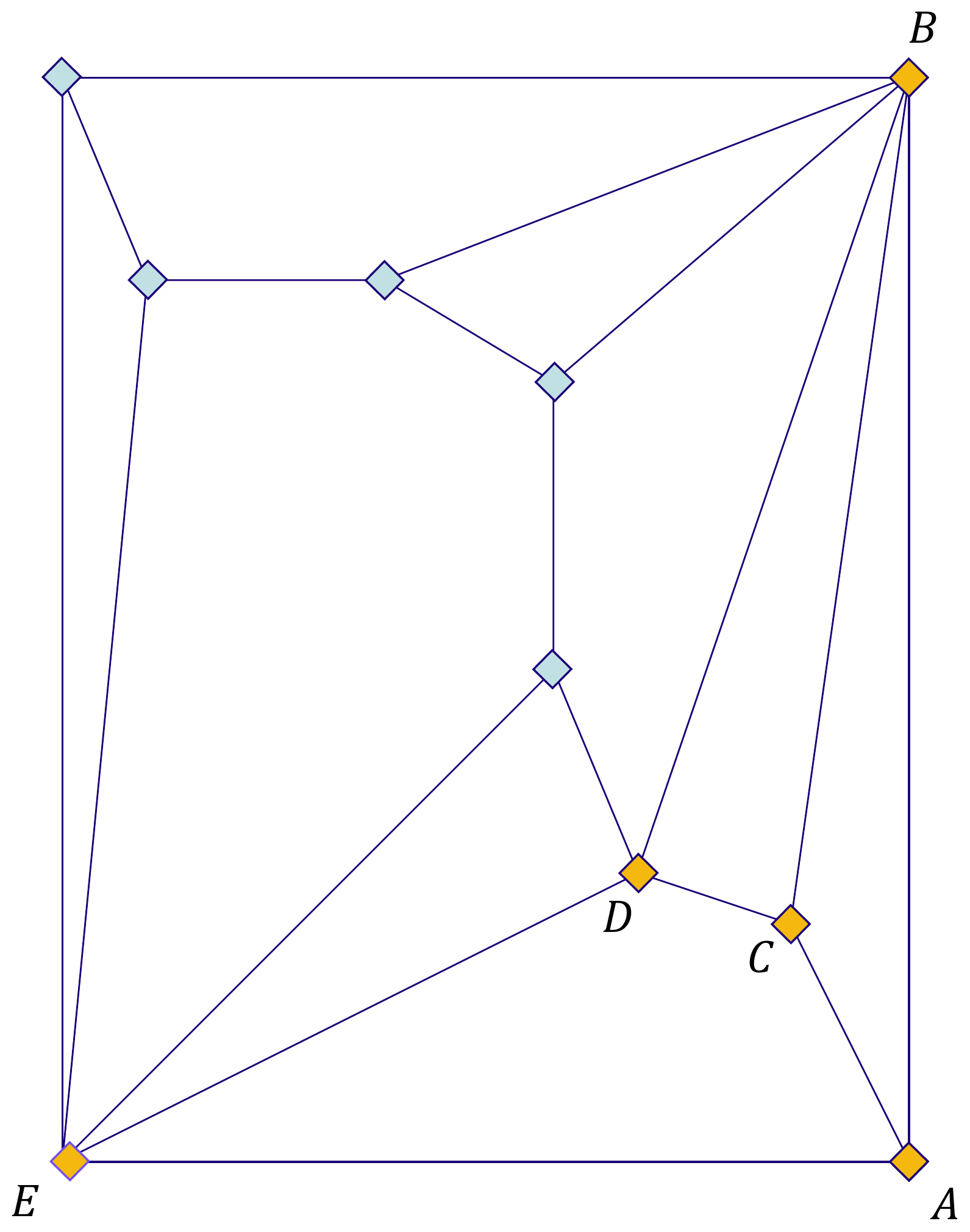}

\caption{\label{rys1} Two indecomposable examples without many triangles}
\end{figure}

To show how powerful Proposition \ref{extend} is, we note that it guarantees indecomposability of any polyhedron whose graph is either of those shown in Figure \ref{rys1}. In neither example is there a chain of triangles touching every face.

For the first example, begin with the edge $AB$, which is indecomposable, then successively add the vertices $C,D, E$ and $F$. Each additional vertex is adjacent to two of the preceding ones, so the resulting geometric graph is indecomposable. Since it touches every face, the polyhedron is indecomposable. The second example is even quicker; beginning with the edge $AB$, it is enough to
 add the vertices $C,D$ then $E$.

A similar argument  also gives a particularly easy proof of the indecomposability of the example in \cite[\S6]{K}. Further applications are given in \cite{BY}.

Kallay \cite[Theorem 8]{K} showed that if two indecomposable graphs have two common vertices, then their union is indecomposable.  A prime example for this result is the 199$^{th}$ polyhedron in the catalogue \cite{BD}, which will be discussed again in the next section. We have drawn it here so that the blue edges are the union of a chain of 3 triangles, and the red edges are the union of another chain of 3 triangles. It is clear that the resulting two indecomposable geometric graphs have two vertices (but no edge) in common, and that their union contains every vertex.

\begin{figure}
\hspace{0.25\textwidth}
\includegraphics[width=0.5\textwidth]{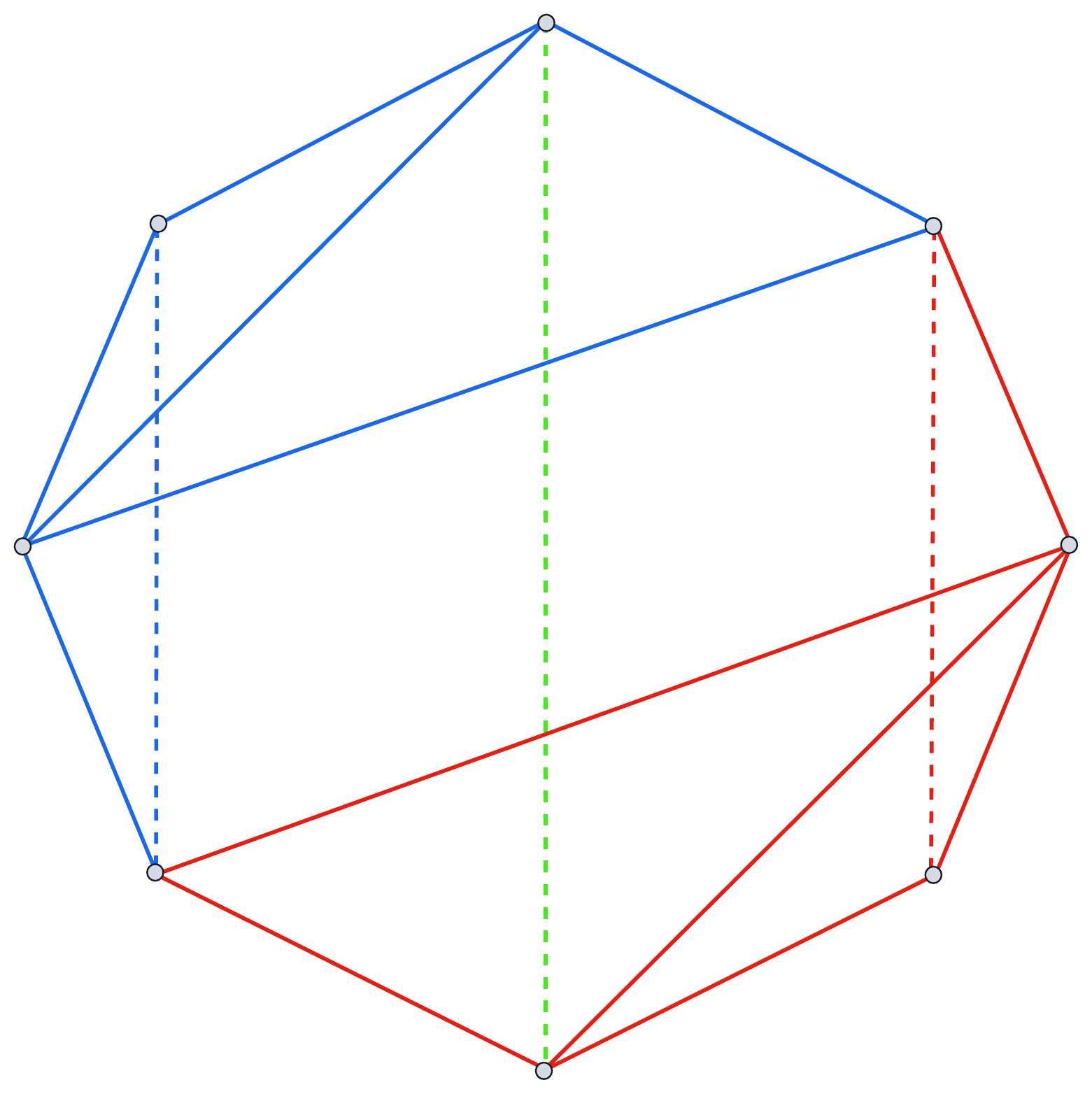}

\caption{\label{rys2} BD199 is indecomposable}
\end{figure}

Our next result is a generalization of both Proposition \ref{extend} and \cite[Theorem 8]{K}. Our proof is no different from Kallay's but, as we shall soon see, our formulation is somewhat more powerful.
It is clear from the definition that adding an edge {but no vertex} to an already indecomposable graph preserves its indecomposability. The point of  part (i) is that, with a little care, we can throw away some edges and still preserve indecomposability.  Part (ii) says that if one edge of an indecomposable graph is replaced by another indecomposable graph, then the new graph is indecomposable.

\begin{thm}\label{omitedge}
{\em(i)}
 Suppose that $G_1=(V_1,E_1)$ and $G_2=(V_2,E_2)$ are two geometric graphs in { the same vector space, and} that $V_{12}=V_1\cap V_2$ contains at least two distinct vertices. Let $E_{12}$ be the collection of those edges of $G_1$, both of whose vertices lie in $V_{12}$. Let $G=(V,E)$ be another geometric graph with vertex set $V=V_1\cup V_2$ and whose edge set $E$  contains
$(E_1\setminus E_{12})\cup E_2$.
If both $G_1$ and $G_2$ are indecomposable, then so is $G$.

\noindent
{\em(ii)} Suppose that $G_1=(V_1,E_1)$ and $G_2=(V_2,E_2)$ are two indecomposable geometric graphs, that $V_1\cap V_2$ contains at least two distinct vertices $u$ and $w$.  Define a new geometric graph $G=(V,E)$ with vertices $V=V_1\cup V_2$ and edges $E=(E_1\setminus\{[u,w]\})\cup E_2$. Then $G$ is also indecomposable.
\end{thm}

\proof (i) Let $f\colon V\to  X$ be a decomposing function, {where $X$ is the ambient vector space}. Since $G_2$ is indecomposable, $f|_{V_2}$ must be the restriction of a similarity, i.e. there are a scalar $\al$ and a vector $x$ such that $f(v)=\al v+x$ for all $v\in V_2$. In particular, $f(u)-f(w)=\al(u-w)$ for all $u,w\in V_{12}$ (even when $[u,w]$ is not an edge of $G_2$). Since $E_1\subseteq E_{12}\cup E$, this implies that $f|_{V_1}$ is also a decomposing function, so by hypothesis must also be the restriction of a similarity. Thus there are a scalar $\be$ and a vector $y$ such that $f(v)=\be v+y$ for all $v\in V_1$.

Now, fix distinct $u,w\in V_{12}$. Consistency requires
$\al u+x=\be u+y$ and $\al w+x=\be w+y$, which quickly forces $x=y$ and $\al=\be$. Thus $f$ is a similarity.

(ii) In the notation of  part (i), we clearly have $(E_1\setminus E_{12})\cup E_2\subseteq E$.
Note that we make no assumption about whether the edge $[u,w]$ belongs to either $G_1$ or $G_2$.\qed

\medskip
Recall that a graph $G$ is called a {\it cycle} if $|V|=k\ge 3$ and $V$ can be ordered as $\{v_1,\ldots , v_k\}$,
so that $E=\{\{v_1, v_2\}, \ldots,\{v_{k-1}, v_k\}, \{v_k,v_1\}\}$. The number $k$ is said to be the {\it length} of the cycle. The next result is a rewording of \cite[Proposition 2]{PY}. Once formulated it is easy to prove, yet surprisingly useful. The 3-dimensional case has already been used in \cite[\S4]{PY}. We state it explicitly here, since we will use both the  3-dimensional and higher dimensional versions in the next sections.

\begin{pr}\label{fourcycle} Any  cycle,  whose vertices are  affinely independent, is an indecomposable geometric graph. In particular, a polytope will be indecomposable, if its skeleton contains a cycle, whose  vertices are not contained in any affine hyperplane, and which touches every facet.
\end{pr}

The next result indicates further how indecomposability of a graph can be established by considering smaller subgraphs.

\begin{thm}\label{newer}

{\em (i)} Let $H=(V,E)$ be an indecomposable geometric graph and  for each $e=[u,v]\in E$, let $G_e=(V_e, E_e)$  be an indecomposable geometric graph containing both vertices $u,v$.  Then the union $\bigcup_e G_e$  is an indecomposable geometric graph.

\noindent
{\em(ii)} Let $G_1=(V_1,E_1),G_2=(V_2,E_2),\ldots,G_n=(V_n,{E_n})$ be indecomposable geometric graphs and  let $v_1,v_2,\ldots,v_n$ be a collection of affinely independent vertices. Set $G_0=G_n$, $v_0=v_n$ and suppose $v_i\in V_i\cap V_{i-1}$ for each $i$. Then the union  $G_1\cup G_2\cup\ldots\cup G_n$ is indecomposable.

\noindent
{\em (iii)} Let $P$ be a polytope, and let $G_1=(V_1,E_1),G_2=(V_2,E_2)$ be two indecomposable subgraphs of the skeleton of $P$, with $V_1\cap V_2\ne\emptyset$, and suppose that $V_1\cup V_2$ contains all but at most $d-2$ vertices of $P$. Then $P$ is indecomposable.

\end{thm}

\proof  (i) Just apply {Theorem} \ref{omitedge}(ii) successively, replacing each edge $e$ of $H$ with the graph $G_e$.

(ii) The graph with vertices $v_1,\ldots , v_n$ and
edges $[v_1, v_2]$, $\ldots$,$[v_{n-1}, v_n]$, $[v_n,v_1]$ is indecomposable by Proposition \ref{fourcycle}. Now we just apply (i).

(iii) If $V_1\cap V_2$ contains two or more elements, the conclusion follows from Theorem \ref{omitedge}(i). So we assume that $V_1\cap V_2$ contains a unique element, say $v_2$. Set $V_i'=V_i\setminus\{v_2\}$, and $C=V\setminus(V_1'\cup V_2')$. Then $C$ contains at most $d-1$ vertices, so their removal from the graph of $P$ will not disconnect it. Since $V_1'$ and $V_2'$ are disjoint, there must then be an edge between them, say between  $v_1\in V_1'$ and $v_3\in V_2'$. Letting $G_3$ be the graph with the single edge $[v_1,v_3]$, we can apply (ii) with $n=3$. (We cannot claim $G_1\cup G_2$ is indecomposable.)
\hfill $\diamondsuit$

\begin{figure}
\hspace{0.06\textwidth}
\includegraphics[width=0.4\textwidth]{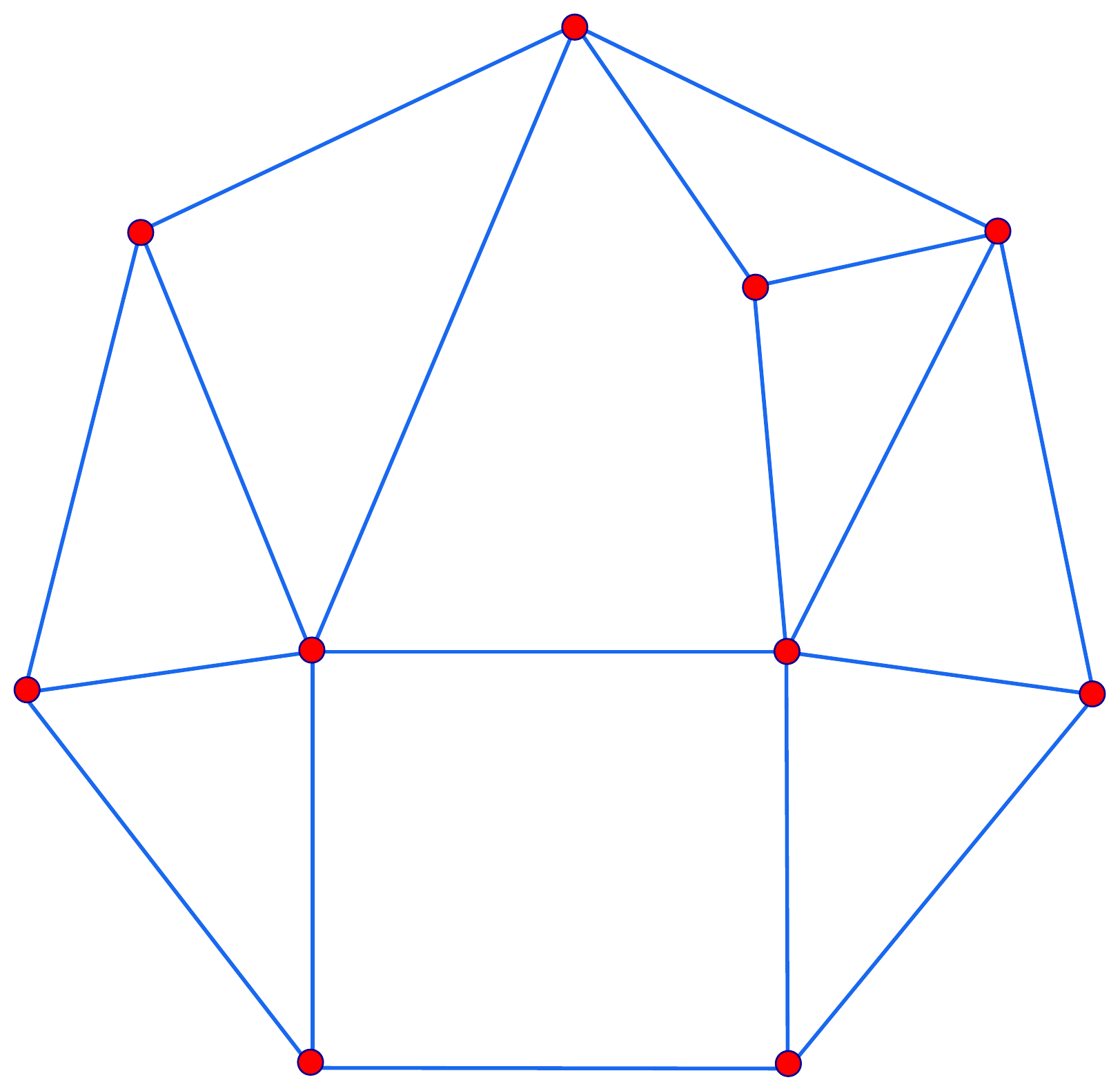}
\hspace{0.08\textwidth}
\includegraphics[width=0.4\textwidth]{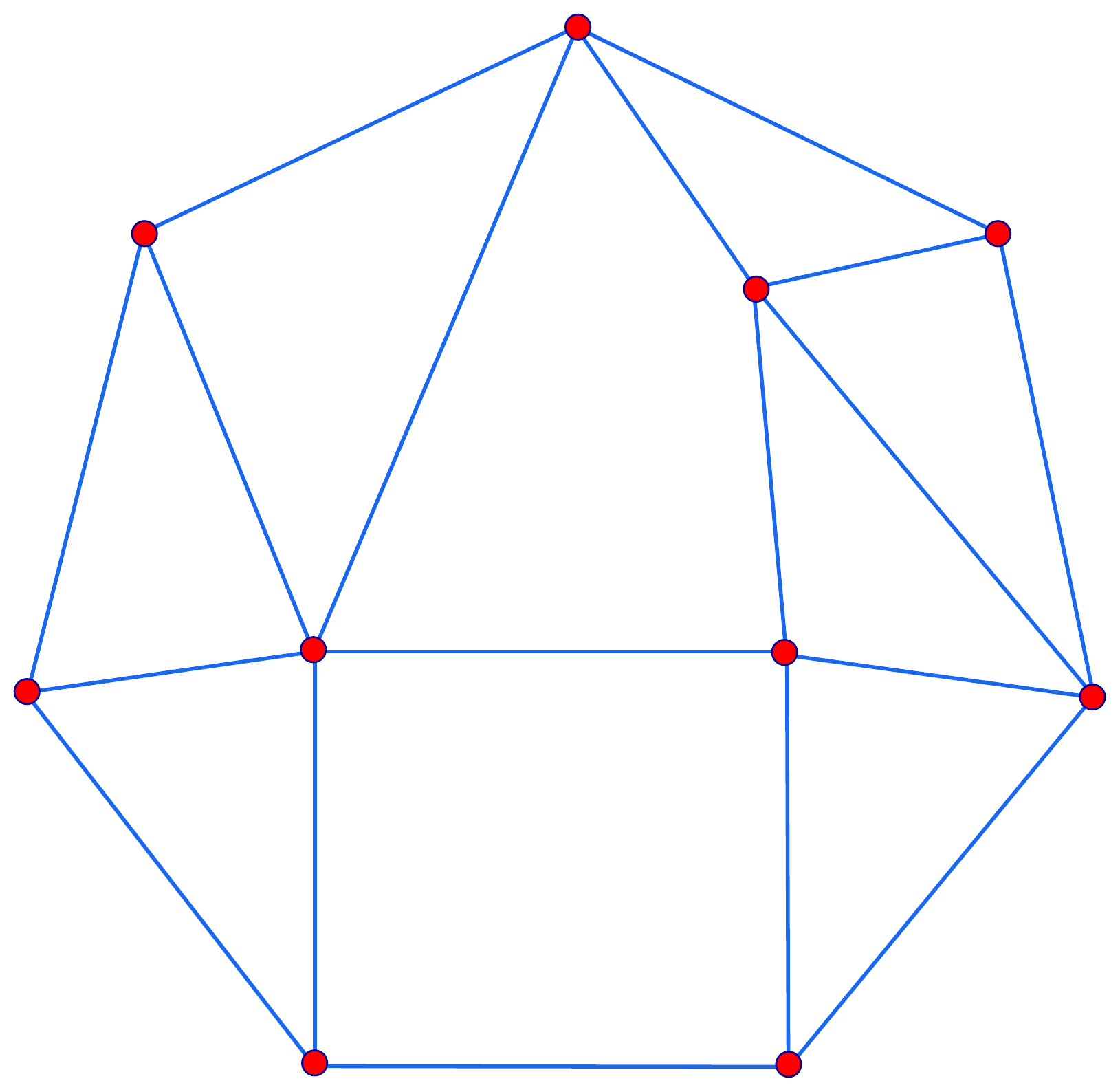}

\caption{\label{rys3} Another two examples }
\end{figure}

\medskip
It is easy to see that \cite[Theorem 9]{K} is precisely the case $n=3$ of part (ii), and that \cite[Theorem 10]{K} is implied by the case $n=4$.  Part (iii) is a strengthening of \cite[Corollary 8.6]{Smi}, where it is assumed that every vertex of $P$ lies in $V_1\cup V_2$. We will indicate the strength of this with another two examples. In both polyhedra whose graphs are shown in Figure \ref{rys3}, we can take $G_1$ and $G_2$ as chains of triangles. Only part (iii) of the preceding Theorem is capable of proving the indecomposability of these two polyhedra.

Sufficient conditions for decomposability  are not so common. The following result was proved without statement by Shephard \cite[Result (15)]{Sh}.  More precisely, he made the stronger assumption that every vertex in $F$ had degree $d$; however, his proof also works in the formulation presented here. It may be interesting to present a proof using decomposing functions.

\begin{pr}\label{shp}
 A polytope $P$ is decomposable whenever there is a facet $F$ such that every vertex in $F$ has a unique neighbor outside $F$, and $P$ has at least two vertices outside $F$.
\end{pr}

\proof
Let $y$ be a support functional for $F$. We may suppose that $y(F)=\{1\}$ and that $y(x)<1$ for all other $x$ in the polytope.  Label the vertices of $F$ as $v_1, \ldots,v_n$. For each vertex $v_i$ of $F$,  denote by $w_i$ the unique vertex which is adjacent to $v_i$ but not in $F$. Set $\al=\max_{i=1}^n y(w_i)$; clearly $\al<1$. For each $i$, let $x_i$ be the unique point on the edge $[v_i,w_i]$ satisfying $y(x_i)=\al$.

Now define a function $f$ by $f(v_i)=x_i$ and $f(v)=v$ for all other vertices. Clearly $f(v)-f(w)=v-w$ whenever both vertices are outside $F$. Since $f(v_i)-f(w_i)=x_i-w_i$ and $x_i$ is a convex combination of $v_i$ and $w_i$, the  condition  for a decomposing function is also satisfied when one vertex lies in $F$. What if both vertices lie in  $F$?   Fix two adjacent vertices $v_i, v_j$  in $F$, and consider a 2-face containing them but not contained in $F$. This face must contain $x_i$ and $x_{j}$. Since
$y(v_i)=y(v_j)\ne y(x_i)=y(x_j)$, the line segments
$[v_i,v_j]$ and $[x_i,x_j]$ must be parallel; we do not claim that $[x_i,x_j]$ is an edge of $P$. Then $f(v_i)-f(v_j)=x_i-x_j$ is a non-negative multiple of $v_i-v_j$. So $f$ is a decomposing function, as anticipated.

Finally, $f$ is not a similarity, because it coincides with the identity function at all the  vertices outside $F$ but is not equal to the identity function. Decomposability follows. \qed

\medskip
The next section requires the 3-dimensional case of the following result. It is not difficult, but appears to be new, so we state it in full generality.

\begin{pr}\label{easy}
Let $F$ be a  facet of a polytope $P$. Suppose that $F$ is indecomposable. Let $Q$ be obtained from $P$ by stacking a pyramid on $F$. Then $P$ is decomposable if and only $Q$ is decomposable.
\end{pr}

\proof Let $V$ be the vertex set of $P$, $u$ the unique vertex of $Q$ not in $P$, $S$ the pyramid being glued onto $F$ and {$X$ the ambient vector space}. It is easy to show that any decomposing function defined on $F$ has a unique extension to $S$.

If $P$ is indecomposable, so is its graph, $G(P)$. Let $f:V\cup\{u\}\to X$ be a decomposing function for $Q$. Then $f|_V$ is a decomposing function for $P$, so we find $\al,x$ so that $f(v)=\al v+x$ for all $v\in V$. By the previous paragraph, $\al u+x$ is the only conceivable value for $f(v)$. Thus $Q$ is indecomposable.

Conversely, suppose $Q$ is indecomposable. Let $g:V\to X$ be any decomposing function for
$P$. Again by the first paragraph, there a unique decomposing function  $f:V\cup\{u\}\to X$ which extends $g$. Since $f$ must be the restriction of a homothety on $X$, so must $g$.  \qed

\section{Polyhedra with 15 edges}

If a given polyhedron has $V$ vertices, $E$ edges and $F$ faces, then Euler's relation $E=V+F-2$ suggests that the number of edges is a reasonable measure of its complexity. Accordingly, we gave in \cite{Y} the complete  classification, in terms of decomposability, of the 58
 combinatorial types of polyhedra with 14  edges. The classification of the 44 types of polyhedra with 13 or fewer edges was essentially known \cite[\S6]{S}.

Indeed, there are only four types of polyhedra with 6, 8 or 10 edges, and they are easily seen to be indecomposable.
No polyhedron can have 7 edges. Besides the triangular prism, the only other polyhedron with 9 edges is the triangular bipyramid, which is obviously indecomposable.

A thorough study of this topic had already been made by Smilansky, who showed \cite[Theorem 6.7]{S} that a polyhedron is decomposable  if there are more vertices than faces; and that a polyhedron is indecomposable if $F\ge2V-6$. As remarked in \cite[p 719]{Y}, simply knowing the values of $F$ and $V$ is then enough to decide decomposability in all cases when $E\le11$ or $E=13$. The examples with 12 edges were discussed in more detail in \cite{Y},  but the results were obviously known to Smilansky. (The only example whose indecomposability is not clear from classical triangle arguments is \cite[figure 2]{Y}, and its indecomposability is guaranteed by \cite[Corollary 8.6]{Smi}.)

The classification of polyhedra with 14 or fewer edges incidentally completed the classification of all  polyhedra with 8 or fewer faces. We should recall that two polytopes are said to be {\it combinatorially equivalent} if their face lattices are isomorphic. In three dimensions, Steinitz's Theorem assures us that two polyhedra are combinatorially equivalent as soon as we know that their graphs are isomorphic. In many cases, two polytopes with the same combinatorial type will either both be decomposable or both be indecomposable. Smilansky \cite[\S 6]{S} first announced that this is not so for polyhedra with 14 edges, and some explicit examples were given in \cite{Y}.

We push this project a bit further in this section by  completing the classification of the 158 combinatorial types of polyhedra with 15 edges. This also  completes the classification of the 301 combinatorial types of polyhedra with 8 or fewer vertices. We also note the indecomposability of all higher dimensional polytopes with 15 or fewer edges.

The aforementioned results of
Smilansky imply that a polyhedron is decomposable  if  $(V,F)$ is either $(10,7)$ or $(9,8)$ and that a polyhedron is indecomposable if  $(V,F)=(7,10)$; these three cases account for 84 combinatorial types. Our assumption that $V+F=17$ then tells us that the only case remaining is $V=8, F=9$.

There are 74 combinatorial types of polyhedra with 8 vertices and 9 faces, which were first described verbally, but not visually, by Kirkman \cite[pp 362--364]{Ki}.
It is possible to use computers to generate diagrams of such polyhedra, but we are dealing with a relatively small number of polyhedra, so it is simpler to use a published catalogue. The only one for this class seems to be that of Britton and Dunitz \cite{BD}. They exhibited diagrams of all the 301 combinatorially distinct types of polyhedra with up to 8 vertices. On their list, those with 8 vertices and 9 faces are numbers 129 to 202 in \cite[Fig. 5]{BD}.

Of these, we will see that most are indecomposable because they have sufficiently many triangles, and
2 are obviously decomposable  thanks to Proposition \ref{shp} (in the simplest geometric realizations, because they have a segment as a summand).
The remaining 6 are also indecomposable but arguments using triangular faces alone don't work; we need to use the results from \S1 to establish their indecomposability.

\begin{thm}\label{main} Of the 74 types of polyhedra with 9 faces and 8 vertices, only 2 types are decomposable, 66 types are indecomposable by classical arguments,  and the remaining 6 require  some results from \S1 to establish their indecomposability. More precisely:

\noindent
{\em(i)} Polyhedra numbers 182 and 198 (on the list of Britton and Dunitz) are decomposable.

\noindent
{\em (ii)} Altogether, 66 are indecomposable by virtue of having a  connected chain of triangular faces. Specifically, we mean those numbered 129--172, 174--178, 180, 181, 183--186, 188, 189, 191, 193--197, 200--202.

\noindent
{\em(iii)} The other 6, namely numbers 173, 179, 187, 190, 192 and 199, are indecomposable thanks to  either Proposition \ref{extend} or \ref{fourcycle} or Theorem \ref{omitedge}.

\end{thm}

\pf\
We begin with the decomposable examples.
As remarked in the opening paragraph, a triangular prism is the Minkowski sum of a triangle and a line segment. If we glue a tetrahedron onto one end of a prism, we obtain the {\it capped prism}, which is decomposable as the sum of a tetrahedron  and a line segment. Better still, Proposition \ref{shp} guarantees that any polyhedron combinatorially equivalent to this  will be decomposable.

There are two ways to glue a second tetrahedron onto the capped prism. Either we glue it onto one face of the first tetrahedron, or we glue it onto the remaining triangular face of the original prism. In both cases, we obtain a polyhedron with 8 vertices and 9 faces. All three are pictured here, the latter two being numbers 182 and 198 respectively from the list of \cite{BD}.

\begin{figure}[!h]
\includegraphics[width=0.3\textwidth]{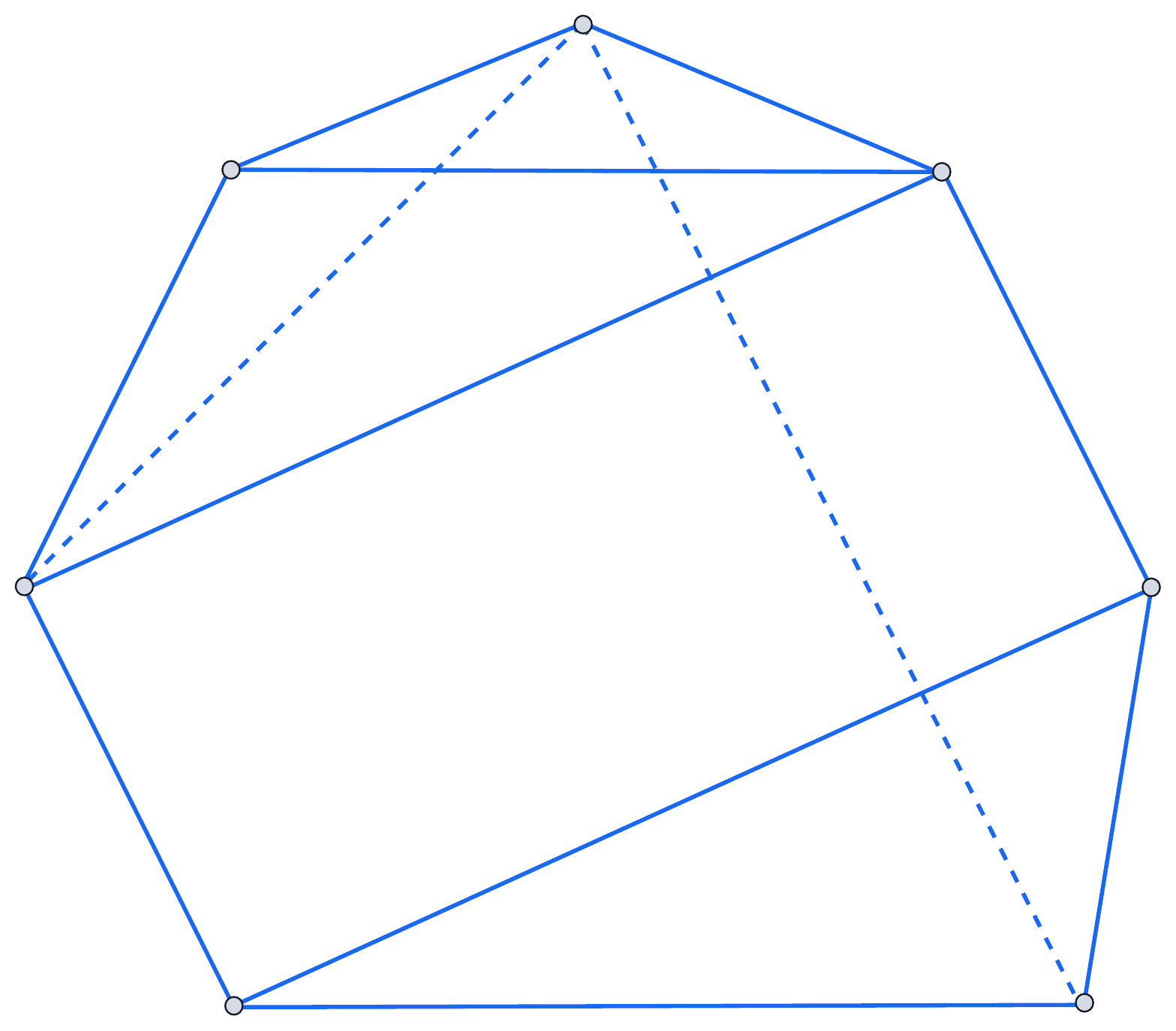}
\hfill
\includegraphics[width=0.3\textwidth]{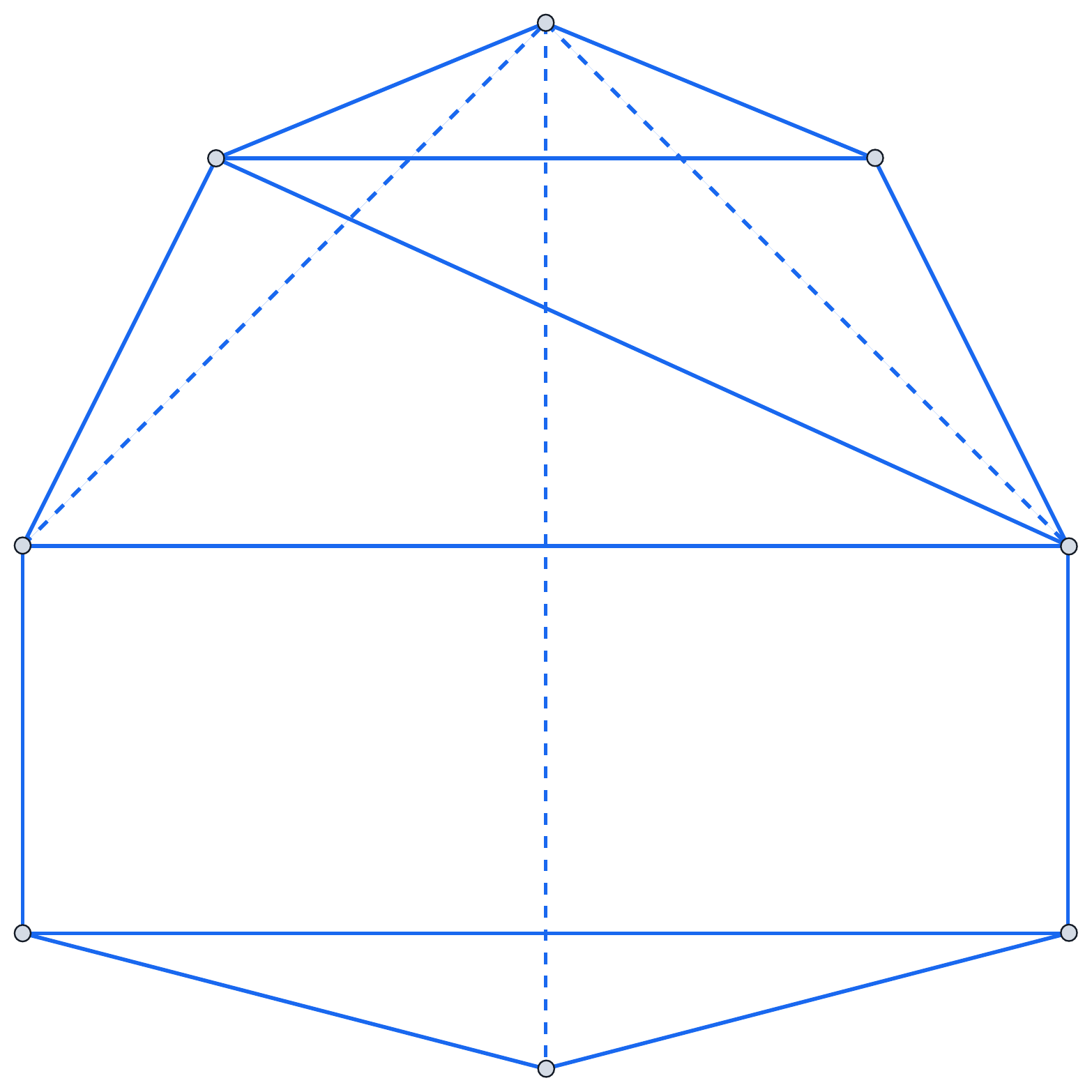}
\hfill
\includegraphics[width=0.3\textwidth]{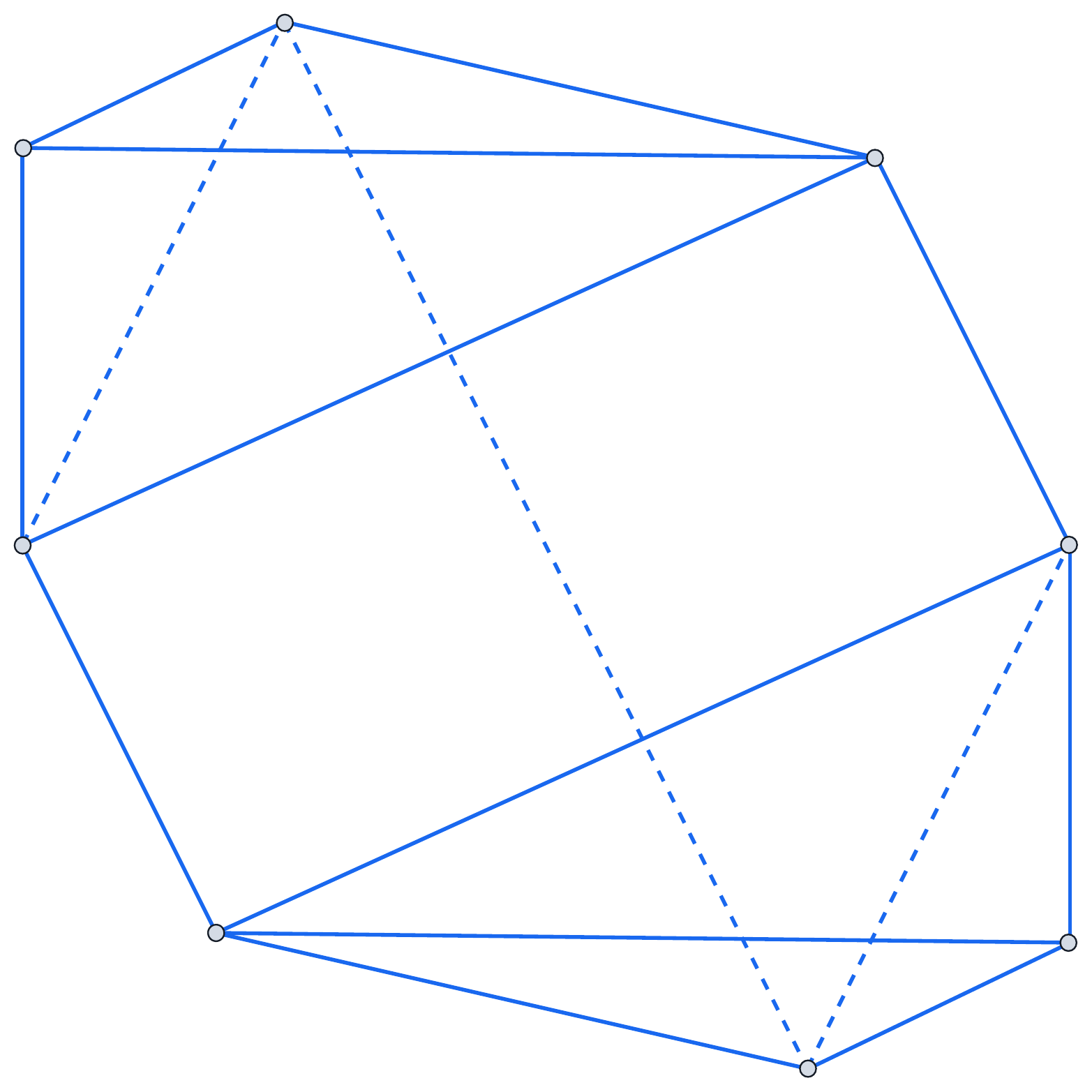}
\caption{The Capped Prism, BD182 and BD198}
\end{figure}

In case the three edges which lie between  pairs of quadrilateral faces are all parallel, each of the latter two polyhedra will be decomposable, being the sum of a triangular bipyramid and a line segment. Since one of them has vertices of degree 5 and the other does not, they are not combinatorially equivalent. This exemplifies the fact that the combinatorial type of two polyhedra does not determine the type of their sum. Our diagrams are not identical to those in \cite{BD}; we have drawn them slightly differently to emphasize their decomposability.

It is true, but not totally obvious, that any polyhedron combinatorially equivalent to these two is decomposable.
Let us prove it.

Proposition \ref{shp} clearly implies that \cite[182]{BD} is decomposable (but not necessarily that a line segment will be a summand).
For \cite[198]{BD}, recall that the capped prism is decomposable, and then apply Proposition \ref{easy}.

Now let us look at the indecomposable examples.
Numbers 129, 130 and 131 each have one hexagonal face and 8 triangular faces. Each of examples 132--155 has one pentagonal face, one quadrilateral face and 7 triangular faces.  Indecomposability of all these examples is assured by our remarks in the previous section, because they have at most two non-triangular faces.

Examples 156--181, 183--197 and 199--202 all have three quadrilateral and six triangular faces. By inspection, all but six
of them (namely 173, 179, 187, 190, 192 and 199) are indecomposable because (some of) their triangular faces can be ordered into a  chain whose union touches every face. We note also that for some examples, the chain of indecomposable triangles does not contain every vertex. (In
particular, 157 and several others each have a vertex which does not lie in any triangular face.) Thus  the weakness of the assumption,
that the chain only touches every face,  is significant.

Proposition \ref{fourcycle} implies the indecomposability of examples
173, 179, 187, 190 and 192 from Britton and Dunitz.  Alternatively, their indecomposability can also be established by Proposition \ref{extend}. None of these examples contains a  connected sequence of triangular faces touching every face, so  some new technique was needed. We present here their
diagrams, with an appropriate 4-cycle highlighted in red. In each case, three vertices of the 4-cycle lie in one face, while the
fourth does not, so the 4-cycle cannot be coplanar. The diagrams make it clear that the 4-cycle touches every face. This time, we
have used the same diagrams as in \cite {BD}, except that for aesthetic reasons we have reversed the front and back faces of 190 and 192.

\begin{figure}
\hspace{0.07\textwidth}
 \includegraphics[width=0.4\textwidth]{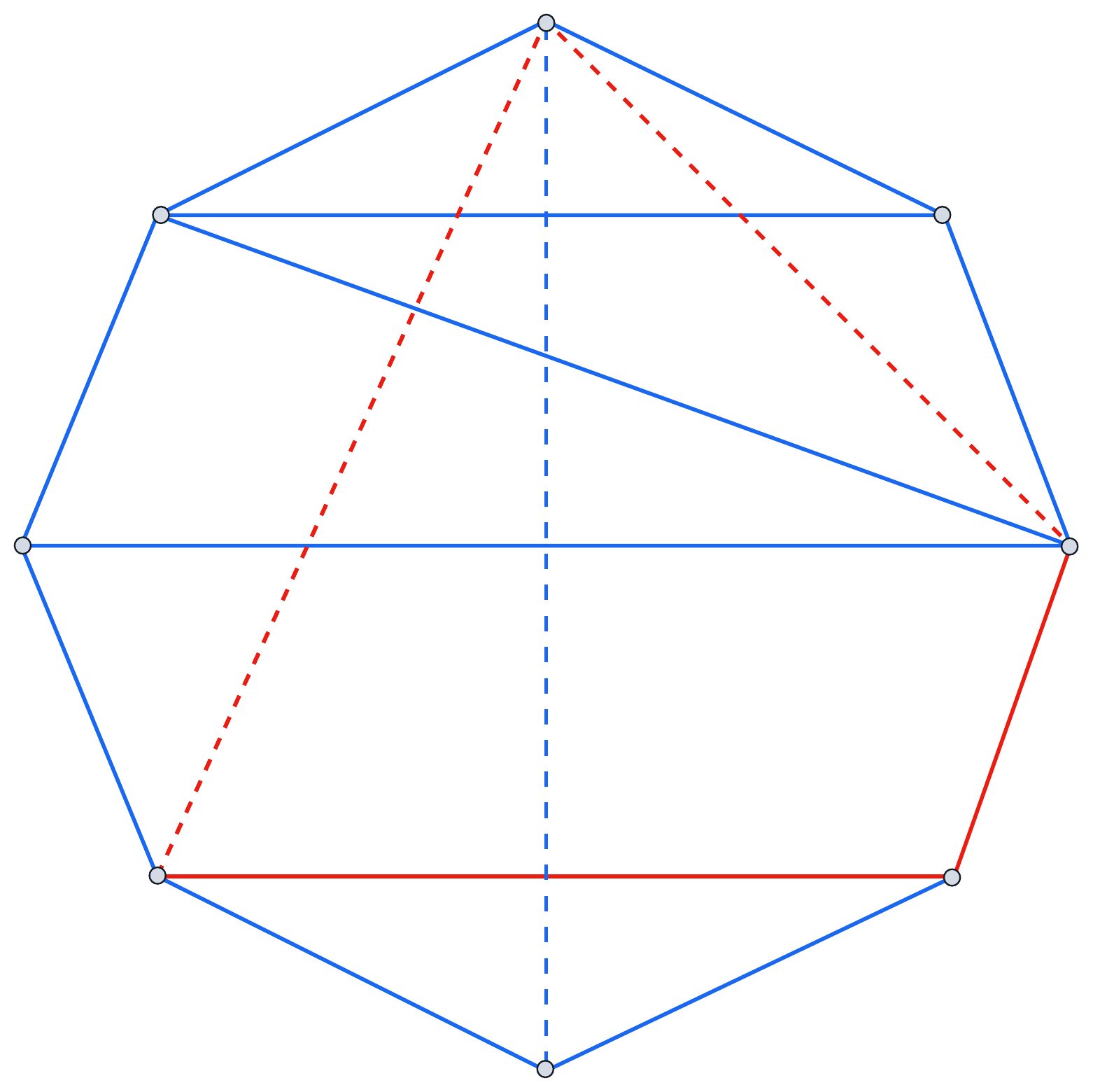}
\hspace{0.06\textwidth}
\includegraphics[width=0.4\textwidth] {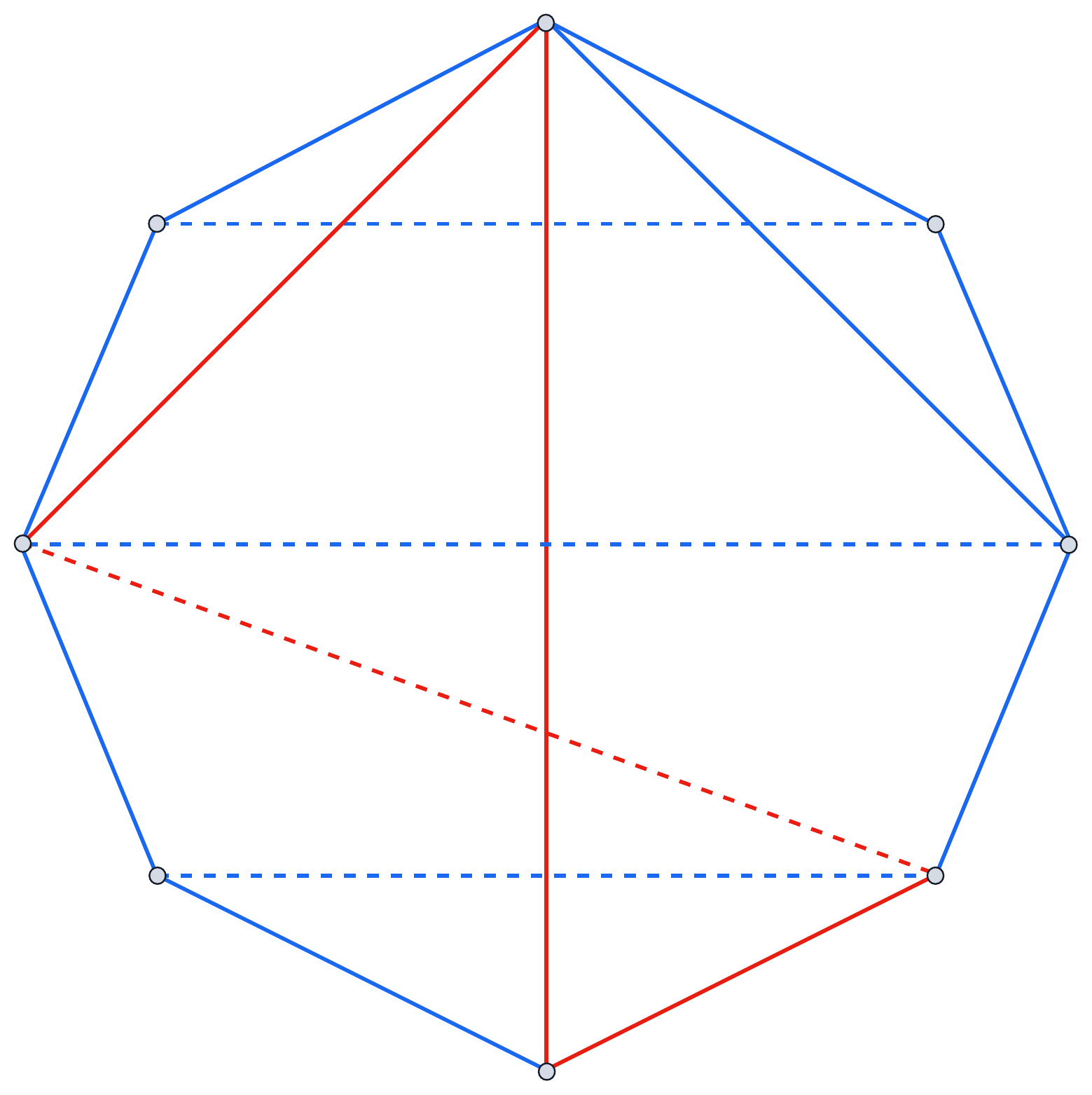}
\caption{BD173 and  BD179}
\end{figure}

\begin{figure}
\includegraphics[width=0.3\textwidth]{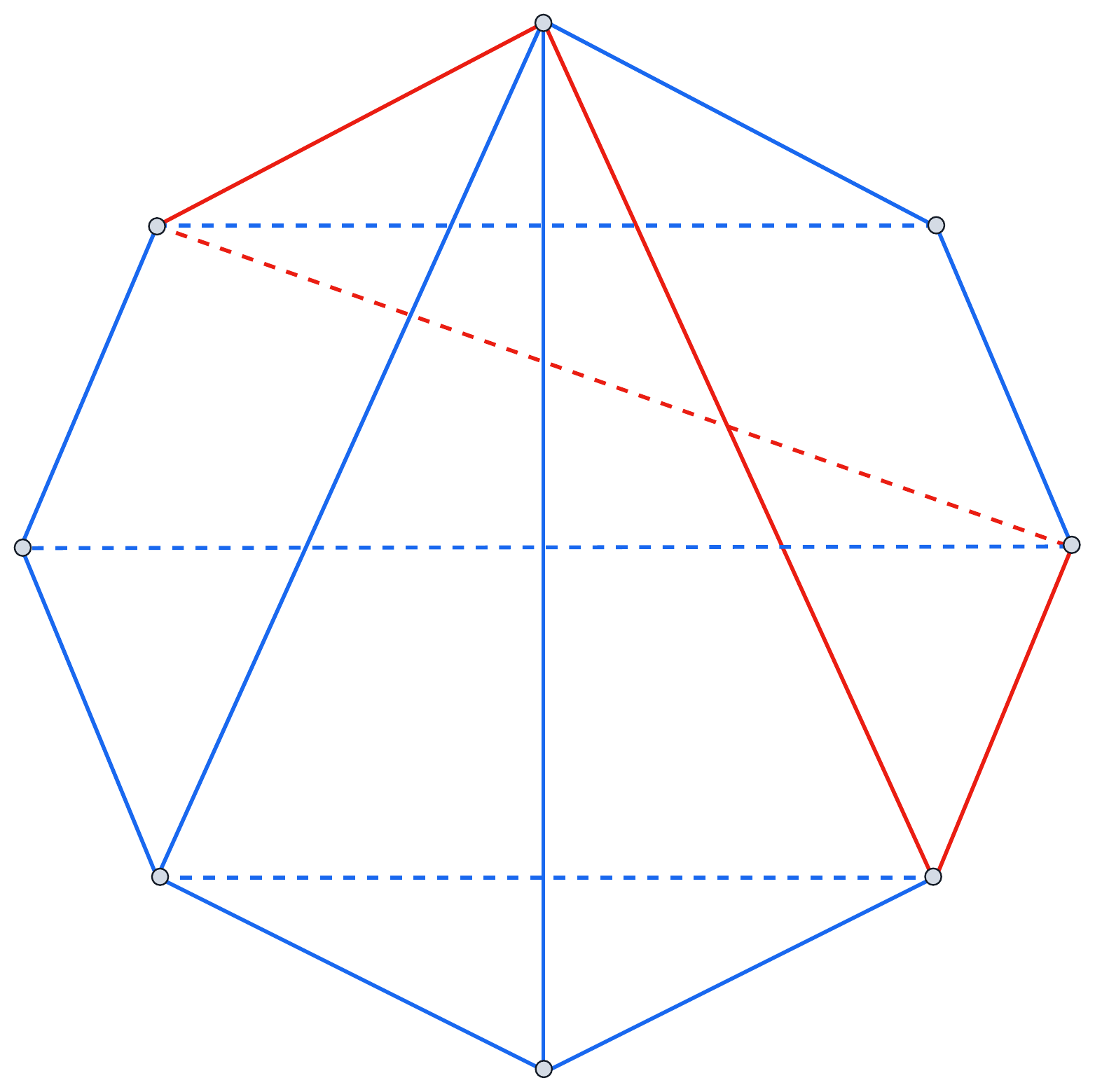}
 \hfill
\includegraphics[width=0.3\textwidth]{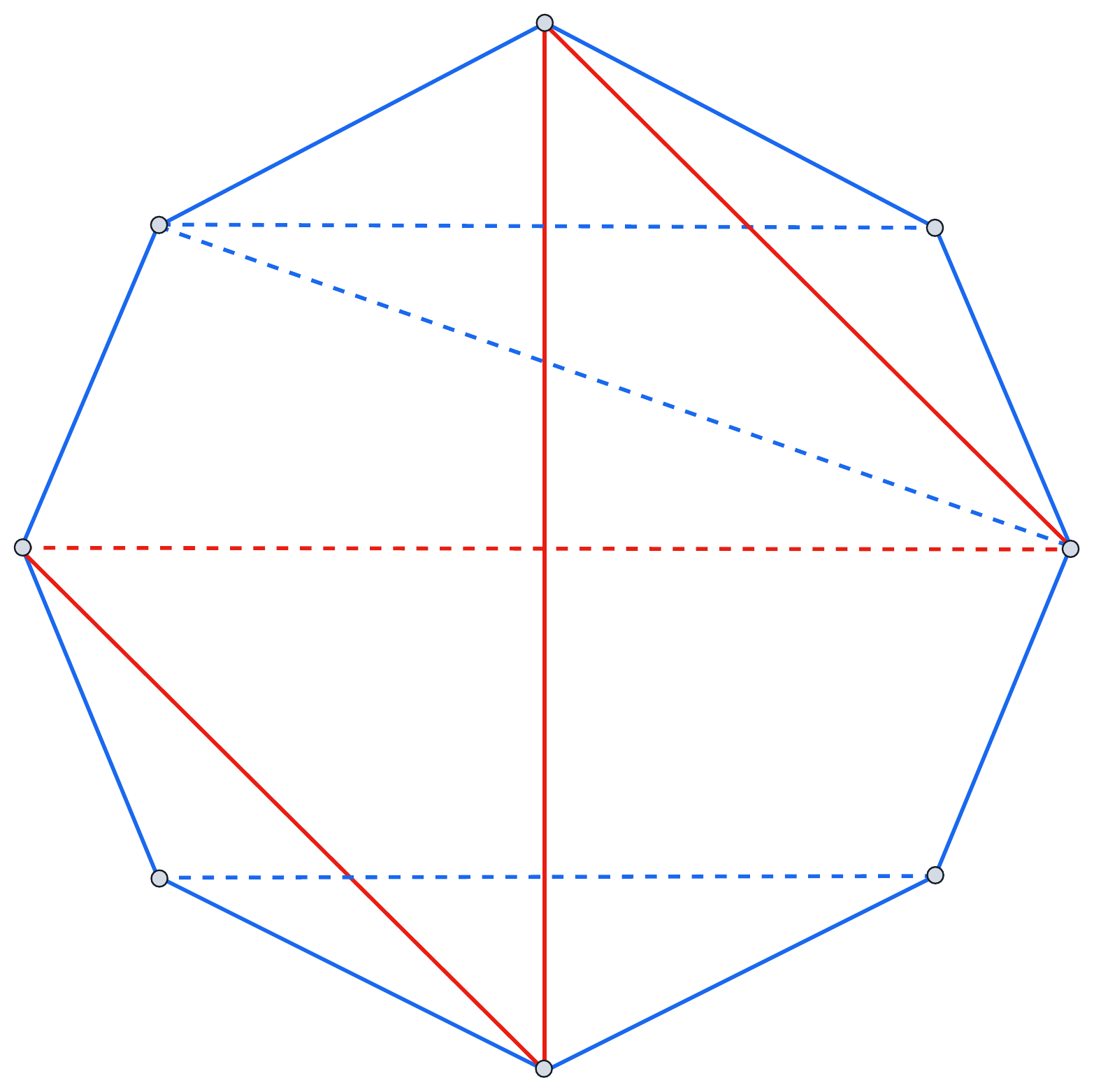}
\hfill
\includegraphics[width=0.3\textwidth]{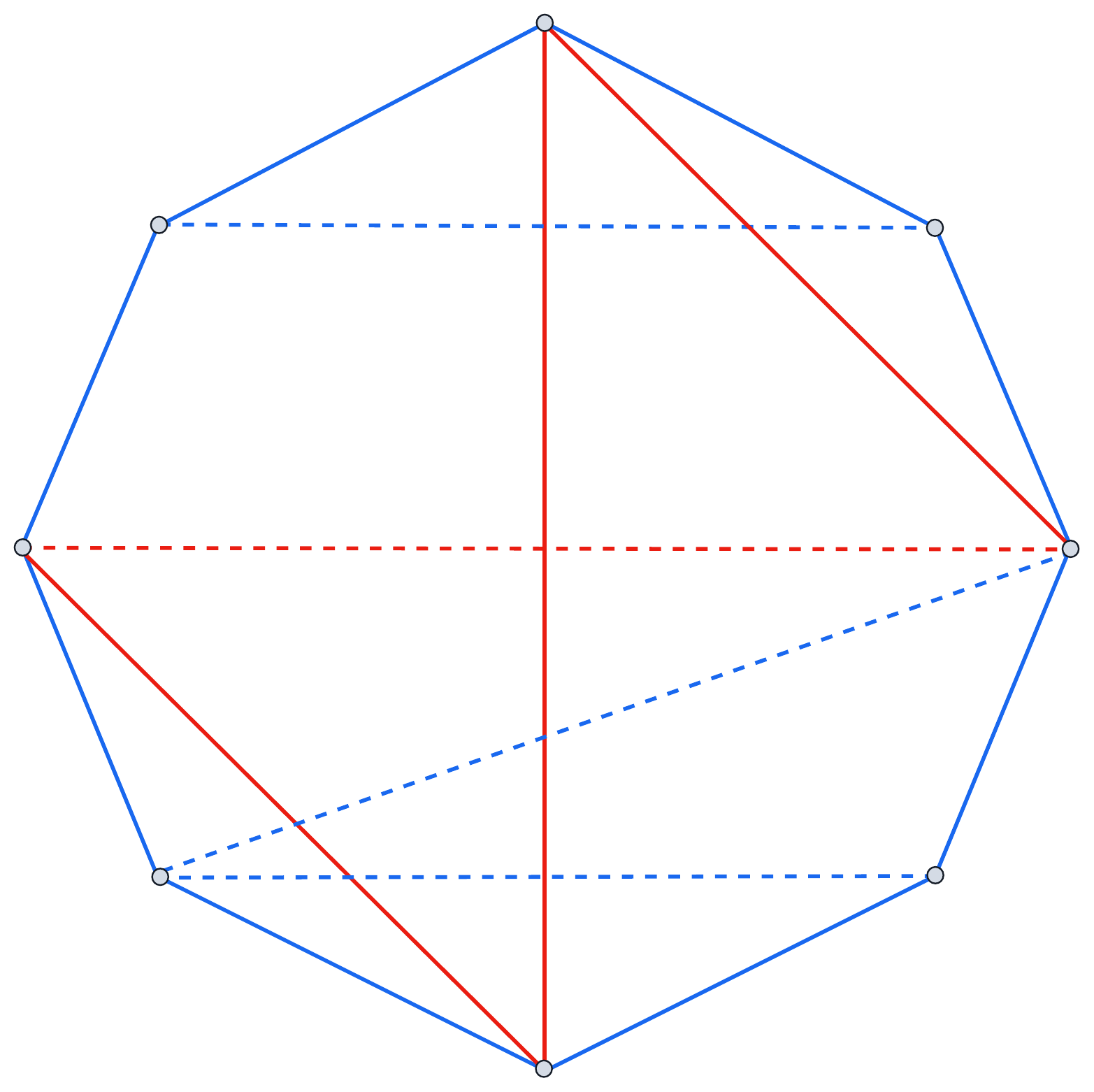}
\caption{BD187, BD190 and BD192}
\end{figure}

Finally, we recall from \S1 that 199 is also indecomposable.  \qed

\medskip
We remark that all higher dimensional polytopes with 15 or fewer edges are indecomposable; this extends \cite[Proposition 2.11]{Y}. The ``smallest" $d$-dimensional polytope, the simplex, obviously has exactly $\h d(d+1)$ edges. So in dimensions 6 and higher, there are in fact no polytopes with 15 or fewer edges. In dimension 5, the only polytope with 15 or fewer edges is the simplex. Two more examples exist in dimension 4, but the next result shows they are both indecomposable.

\begin{pr}\label{fourdim}
The assertion ``every $4$-dimensional polytope with $n$ edges is indecomposable" is true if and only if $n\le15$ or $n=17$.
\end{pr}

\pf Any polytope satisfying these restrictions on $n$ will have at most 7 vertices \cite[10.4.2]{Gr}. This condition forces indecomposability by \cite[Proposition 6]{Yo}.

For the converse, we need to consider various possible values for $n$. We will simply describe the examples, and not verify all the details.

A particularly simple decomposable polytope with  18 edges is the sum of two triangles lying in orthogonal planes.

The sum of a 4-dimensional simplex with a line segment, which is parallel to one 2-face but not parallel to any edge of the simplex, will have 19 edges.

The sum of a 4-dimensional simplex with a line segment, which is not parallel to any proper face of the simplex, will have 20 edges.

Denoting by $e_i$ the usual basis vectors, let $P$ be the convex hull of $\{0,e_1,$  $e_2,e_3,e_4,e_3+e_4\}$. Then the sum of $P$ with the segment $[0,e_1]$ has 22 edges.

The sum of the cyclic polytope $C(6,4)$ with a line segment which is parallel to one of its edges, will have 25 edges.

The sum of a  4-dimensional simplex with a triangle, which is parallel to one of its 2-faces but has the opposite orientation, will have 27 edges.

If $P$ is a polyhedron with $E$ edges and $V$ vertices, then the sum of a $P$ with a line segment (not parallel to the affine hull of $P$) is easily seen to have $2E+V$ edges. (This is equally true in higher dimensions.) The possible values of $E$ and $V$ for polyhedra are well known \cite[\S10.3]{Gr}, and the corresponding values of $2E+V$ account for all remaining values of $n$. \qed

\medskip
In particular, the sum of a tetrahedron with a line segment  has  16 edges. In the next section, we will see that this is (up to combinatorial equivalence) the only example with 16 edges.

\section{Polytopes with not too many edges}

A simplicial prism, i.e. the sum of a segment with a $(d-1)$-dimensional simplex, has $2d$ vertices, $d^2$ edges and $d+2$ facets. These numbers turn out be the minimum possible, for a $d$-dimensional decomposable polytope. In the case of vertices or edges, the prism is (up to combinatorial equivalence) the unique minimiser.

In $d$ dimensions, any polytope has at least $d+1$ facets, and only the simplex has $d+1$ facets. So no non-trivial bound on the number of facets will imply indecomposability. Nor can uniqueness be expected;  a $(d-2)$-fold pyramid over a quadrilateral also has $d+2$ facets. For further examples, see Lemma \ref{prismoid} below.

The conclusions regarding the numbers of vertices and edges  are  more interesting; for edges, this extends Proposition \ref{fourdim}
to higher dimensions.
Proposition \ref{fourcycle} is an essential tool for these. So also  is Gale's result \cite[(14)]{Sh} that any pyramid, i.e. the convex hull of a maximal face and a single point, is indecomposable. This is clear, because every 2-face outside the base must be triangular.

As noted in \cite[Proposition 6]{Yo}, a $d$-dimensional polytope with strictly fewer than $2d$ vertices is automatically indecomposable, and this estimate is the best possible.

We will prove  now that the simplicial prism is the only decomposable $d$-polytope with $2d$ or fewer vertices, before the corresponding result about edges. We have learnt recently that this result was first proved by Kallay  \cite[Theorem 7.1, page 39]{Kal} but never published; his argument is different, using Balinski's Theorem.

Recall that a $d$-polytope $P$ is simple  if every vertex is simple, i.e. has degree $d$. Clearly every simple $d$-polytope, other than a simplex, is decomposable.

\begin{thm}\label{kallaythesis}
 Let $P$ be a decomposable $d$-dimensional polytope with $2d$ or fewer vertices. Then $P$ is combinatorially equivalent to the sum of a line segment and a $(d-1)$-dimensional simplex (and hence has precisely $d^2$ edges).
\end{thm}

\pf  The 2-dimensional case is almost obvious and the 3-dimensional case is quite easy, from \S2. We proceed by induction on $d$.

So let $P$ be a decomposable $(d+1)$-dimensional polytope with $2(d+1)$ or fewer vertices.

Then some $d$-dimensional facet, say $F$, must be decomposable. Since $P$ is not a pyramid, there must be (at least) two vertices of $P$ outside $F$; this implies that $F$ has at most $2d$ vertices. By the inductive hypothesis, $F$ is combinatorially equivalent to the sum of a line segment and a $(d-1)$-simplex.

This means that $F$ has two faces which are simplices, whose vertex sets $\{v_1,v_2,\ldots,v_d\}$ and $\{w_1,w_2,\ldots,w_d\}$ can be labelled in such a way that $v_i$ is adjacent to $w_i$ for each $i$. In particular $F$ has $2d$ vertices and $d^2$ edges.

Furthermore there must be
precisely two vertices of $P$ outside $F$, say $x$ and $y$.

Suppose that one of them is adjacent  to vertices in both simplices, say $[x,v_i]$
and $[x,w_j]$ are both edges of $P$ for some $i$ and $j$. A routine degree argument shows that $x$ is adjacent to at least two vertices in one simplex, so without loss of generality $i\ne j$. We may renumber the vertices so that $i=1$,  $j=d$. But then
$$\{v_1,v_2,\ldots, v_d,w_d,x\}$$
will be an affinely independent $(d+2)$-cycle. It
touches every facet, since  $P$ has only $2d+2$ vertices. This contradicts our assumption that $P$ is decomposable.

Thus each of $x,y$ is adjacent to vertices in only one simplex, say $x$ is not adjacent to any $w_j$ and $y$ is not adjacent to any $v_i$. Since all vertices have degree at least $d+1$,
it follows that $x$ is adjacent to each $v_i$, $y$ is adjacent to each $w_j$, and $x$ and $y$ are adjacent to each other. This means that the skeleton of $P$ is isomorphic to the skeleton of the sum of a line segment and a simplex.

Now observe that $P$ is simple and so is in fact combinatorially equivalent to the sum of a line segment and a simplex, thanks to a result of Blind and Mani \cite{BML}.
\qed

\medskip
Proposition \ref{shp} implies that if we cut any vertex from any polytope, the resulting polytope will be decomposable. This makes it easy to construct decomposable polytopes with any number of vertices greater than $2d$. On the other  hand, Proposition \ref{fourdim} asserts that there are gaps in the possibe numbers of edges of decomposable polytopes, at least in dimension 4. We show now that this is also true in higher dimensions. In fact, a decomposable $d$-dimensional polytope with strictly less than $d^2+\h d$  edges must be combinatorially equivalent to a prism; this is  an easy consequence of Theorem \ref{kallaythesis}. With some additional material, we can prove a stronger result.

We will first examine the existence of simple polytopes with less than $3d$ vertices. Being decomposable,  Theorem \ref{kallaythesis}  implies that no simple $d$-polytope has between $d+1$ and $2d$ vertices. This also follows from Barnette's Lower Bound theorem.  For results concerning higher numbers of vertices, see \cite{P} and the references therein.

We denote by $\Delta_{m,n}$ the sum of an $m$-dimensional simplex and an $n$-dimen\-sio\-nal simplex lying in complementary subspaces. It is routine to check that $\Delta_{m,n}$ is a simple $(m+n)$-dimensional polytope with $(m+1)(n+1)$ vertices, $\h(m+n)(m+1)(n+1)$ edges and $m+n+2$ facets. We denote by $W_d$ the result of cutting a vertex from a $d$-dimensional simplicial prism $\Delta_{1,d-1}$. This simple polytope has $3d-1$ vertices, $\h d(3d-1)$ edges, and $d+3$ facets, comprising 2 simplices, 2 prisms and $d-1$ copies of $W_{d-1}$. In dimension 3, $W_3$ is simply the 5-wedge.

\begin{lem}\label{prismoid}
(i) The (combinatorial types of) simple $d$-dimensional polytopes with $d+2$ facets are precisely  the polytopes $\Delta_{k,d-k}$ for $1\le k\le\h d$.

(ii) Up to combinatorial equivalence, the only simple $d$-dimensional polytopes with fewer than $3d$ vertices are the simplex $\Delta_{0,d}$, the simplicial prism $\Delta_{1,d-1}$, the polytope $\Delta_{2,d-2}$,  the 6-dimensional polytope $\Delta_{3,3}$, the polytope $W_d$, the 3-dimensional cube $\Delta_{1,1,1}$ and the 7-dimensional polytope $\Delta_{3,4}$.

(iii) For every $d\ne6$, the smallest vertex counts of simple $d$-polytopes are $d+1$, $2d$, $3d-3$ and $3d-1$. In dimension 6 only, there is also a simple polytope with $3d-2$ vertices.

\end{lem}

\pf (i) The simplicial polytopes with $d+2$ vertices are described in detail by Gr\"unbaum \cite[\S6.1]{Gr}, and these are their duals.

(ii) Obviously the simplex is the only polytope with $d+1$ (or fewer) facets. Barnette, \cite{B} or \cite[\S19]{Br}, showed that a polytope with $d+4$ or more facets has at least $4d-2\ge3d$ vertices. He also showed that a polytope with  $d+3$  facets has at least $3d-1$ vertices, and that if $d>3$ the only such example with precisely $3d-1$ vertices arises from truncating a vertex from a simplicial prism, i.e. it is $W_d$. If $d=3$, the cube $\Delta_{1,1,1}$ is the unique other example.

We are left with the case of $d+2$ facets. Clearly $\Delta_{1,d-1}$ and $\Delta_{2,d-2}$ have respectively $2d$ and $3d-3$ vertices.

If $3\le k\le\h d$, then $d\ge6$. If $d\ge8$, then $\Delta_{k,d-k}$ has at least $(3+1)(d-3+1)>3d-1$ vertices. If $d=7$, we have the example $\Delta_{3,4}$, which has $20=3d-1$ vertices. If $d=6$, we must also consider $\Delta_{3,3}$, which has $16=3d-2$ vertices.

(iii) This follows immediately from (ii). \qed

\begin{thm}\label{lastest}
 Let $P$ be a decomposable $d$-dimensional polytope with no more than $d^2+\h d$  edges. Then either $P$ is combinatorially equivalent to a simplicial prism $\Delta_{1,d-1}$ (and hence has precisely $d^2$ edges), or $d=4$ and $P$ is combinatorially equivalent to $\Delta_{2,2}$.
\end{thm}

\pf  A $d$-dimensional polytope with $2d+1$ or more vertices must have at least $\h(2d+1)d$ edges.

So if $P$ has $2d+1$ vertices, it must be simple, and Lemma \ref{prismoid} implies that $2d+1\ge3d-3$. Thus $d=4$ and $P$ is  $\Delta_{2,2}$.

Otherwise, $P$ has at most $2d$ vertices and the conclusion follows from Theorem \ref{kallaythesis}. \qed

\medskip
In particular, a polychoron with 17 edges is necessarily indecomposable. Gr\"unbaum \cite[p 193]{Gr} showed that there is no polychoron at all with 8 vertices and 17 edges. We finish by using the preceding results to show that this is not an isolated curiosity: in fact, there is no  $d$-dimensional polytope with $2d$ vertices and $d^2+1$ edges for any higher value of $d$. (There are two easy examples when $d=3$; see \cite[Fig. 3]{BD}.)

\begin{lem}\label{technical}
The polytope $\Delta_{2,d-3}$ cannot be a facet of any decomposable $d$-di\-men\-sio\-nal polytope with $3d-4$ vertices.
\end{lem}

\pf We can realize $\Delta_{2,d-3}$ as the convex hull of three $(d-3)$-simplices, say $S, T, U$, all translates of one another, so that the convex hull of any two of them is a facet therein, combinatorially equivalent to  $\Delta_{1,d-3}$. Moreover in each such facet, e.g. $\co(S,T)$, each of the $d-2$ edges joining $S$ and $T$ also belongs to a triangular face whose third vertex lies in $U$.

Suppose that this copy of $\Delta_{2,d-3}$ is a facet of a decomposable polytope $P$ with $3d-4$ vertices. Denote $v,w$ the two vertices of $P$ lying outside this facet. Then $\co(S,T)$ is a ridge in $P$; denote by $F$ the other facet containing it. Then $F$ contains at least one of $v, w$.

In particular, $F$ omits at most $d-1$ vertices of $P$. These $d-1$ vertices cannot form a facet, so $F$ touches every facet. Decomposability of $P$ then implies that  $F$ is also decomposable.

Since $F$ has at most $2d-2$ vertices, it can only be a copy of the prism $\Delta_{1,d-2}$, with one of $v,w$ adjacent to every vertex in $S$ and no vertex in $T$, and the other adjacent to every vertex in $T$ and no vertex in $S$. The same argument applied to $\co(T,U)$ and $\co(S,U)$ quickly yields a contradiction.
\qed

\begin{thm}\label{application}
 Let $P$ be a
  $d$-dimensional polytope with $2d$ vertices and $d^2+1$ or fewer edges. Then either $P$ is combinatorially equivalent to the prism $\Delta_{1,d-1}$ (and hence has precisely $d^2$ edges), or $d=3$.
\end{thm}

\pf If $d$ is 1 or 2, the conclusion is obvious. In case we can establish decomposability, the conclusion will follow from Theorem \ref{lastest}.

If $P$ has exactly $d^2$ edges, then it is simple, hence decomposable by Shephard's result, Proposition \ref{shp}. Since every vertex has degree at least $d$, $P$ cannot have fewer than $d^2$ edges.

We are forced to contemplate the possibility that $P$ has precisely $d^2+1$ edges. Then $P$ is indecomposable by Theorem \ref{lastest}. Since $2E-dV=2$,  there are at most two vertices which are not simple.

Now suppose that  some vertex has degree $d+2$, and choose a facet $F$ not containing $v$. Then $P$ must be a pyramid over $F$, otherwise it would be decomposable by Shephard's result. Then $F$ has $v=2d-1=2(d-1)+1$ vertices, and hence at least $\h (d-1)v=(d-1)^2+\h(d-1)$ edges. Hence $P$ will have at least $(d-1)^2+\h(d-1)+ (2d-1)=d^2+\h d-\h$ edges. The hypothesis then implies that $\h d-\h\le1$. We conclude that $d=3$ and $P$ is a pentagonal prism.

Next consider the case that one vertex $v$ has degree $d+1$ and that all its neighbors are simple vertices. If we cut this vertex from $P$, the resulting facet will be simple and contain $d+1$ vertices. This facet cannot be a simplex, so  Lemma \ref{prismoid} implies that $d+1\ge2(d-1)$, i.e. $d\le3$.

Finally consider the case that $P$ has two adjacent vertices of degree $d+1$. We can find a hyperplane which has this edge on one side, and all other vertices of $P$ on the other side. This divides $P$ into two polytopes, say $Q$ and $R$ respectively, with a common facet $F$. All other vertices are simple, so $F$ will be simple and contain $2d$ vertices. Now $F$ cannot be a simplex or a prism, because it has more than $(d-1)+1$ or $2(d-1)$ vertices. Lemma \ref{prismoid}(iii) then forces $2d\ge3(d-1)-3$, i.e. $d\le6$.

If $d=6$, then $F$ has 12 vertices, and can only be $\Delta_{2,3}$. But $Q$ has 14 vertices, which is impossible according to  Lemma \ref{technical}. If $d=5$, then $F$ is simple and has  $10=3(d-1)-2$ vertices, which according to Lemma \ref{prismoid} is impossible unless $d-1=6$.
Gr\"unbaum \cite[p 193]{Gr} showed that the case $d=4$ is impossible. The only remaining possibility is that $d=3$ and $P$ is combinatorially equivalent to the second last example in \cite[Fig. 3]{BD}.\qed

\medskip

\noindent
{\bf Acknowledgements.} We thank Eran Nevo for assistance in translating reference \cite{Kal}, and Vladimir Fonf for assistance in translating  reference \cite{S}. The second author records his thanks to the University of Zielona G{\'o}ra, for hospitality during his visits in 2008 and 2011.

Wydzia{\l}  Matematyki, Informatyki i Ekonometrii,

Uniwersytet Zielonog\'orski,

ul. prof. Z. Szafrana 4a,

65-516 Zielona G\'ora,

POLAND

e-mail: K.Przeslawski@wmie.uz.zgora.pl
\vskip0.5cm
and
\vskip0.5cm
Wydzia{\l} Matematyki, Informatyki i Architektury Krajobrazu,

Katolicki Uniwersytet Lubelski,

ul.  Konstantyn\'ow 1 H,

20-708 Lublin,

POLAND

\vskip1.3cm

Centre for Informatics and Applied Optimization,

Faculty of Science and Technology,

Federation University,

PO Box 663,

Ballarat, Vic. 3353

AUSTRALIA

e-mail: d.yost@federation.edu.au

\bigskip

\end{document}